\documentclass[11pt]{article}
\usepackage{amsmath,amssymb,amsthm}
\usepackage[margin=1in]{geometry}
\usepackage{graphicx}
\usepackage{longtable}
\usepackage{array}
\usepackage[table]{xcolor}
\definecolor{lvH}{HTML}{1B6F8C}
\colorlet{headbg}{lvH!20}
\colorlet{newbg}{lvH!8}
\colorlet{knownbg}{black!4}
\colorlet{relbg}{white}
\colorlet{conbg}{lvH!8}
\colorlet{stripebg}{lvH!4}
\arrayrulecolor{black!35}
\newcommand{\semrow}[1]{\rowcolor{#1}}
\usepackage{enumitem}
\usepackage{hyperref}

\newtheorem{theorem}{Theorem}
\newtheorem{proposition}[theorem]{Proposition}
\newtheorem{lemma}[theorem]{Lemma}
\newtheorem{corollary}[theorem]{Corollary}
\theoremstyle{definition}

\newtheorem{remark}[theorem]{Remark}

\newcommand{\SOS}{\operatorname{SOS}}
\newcommand{\BSR}{\operatorname{BSR}}

\title{On $z_{2M}(21,4)$: the Hole Extensions of the\\ Chen--Chen Configuration Are Reducible,\\ and the Extremal Configurations of $z_2(21,4)$ Are Not Unique}
\author{
    Carles Mar\'in\\
    {\small Independent researcher}\\
    {\small\tt karlesmarin@gmail.com}
    \and
    Liqun Qi\\
    {\small School of Science, Zhejiang University of Science and Technology, Hangzhou, China}\\
    {\small Department of Applied Mathematics, The Hong Kong Polytechnic University, Hong Kong}\\
    {\small\tt maqilq@polyu.edu.hk}
    \and
    Gaohang Yu\\
    {\small School of Science, Zhejiang University of Science and Technology, Hangzhou, China}\\
    {\small\tt maghyu@163.com}
}
\date{\today}

\begin{document}
\maketitle

\begin{abstract}
The second order extended Zarankiewicz number $z_{2M}(m,n)$ of Qi and Yu
allows a secondary half of a $2$-edge to be shared, while the $1$-edge
graph stays $C_4$-free and extremal and $1$-edges and main halves stay
exclusive.  Whether $z_{2M}(21,4)=55$ is open.  Since deleting a $2$-edge
removes exactly one displayed square, $z_{2M}(21,4)=55$ if and only if
there is no irreducible admissible configuration with exactly $56$
displayed squares, and by the cell count of Qi and Yu and the heredity of
irreducibility such a configuration without empty cells is an extremal
$55$-square configuration of $z_2(21,4)$ plus one $2$-edge at its hole.
For the extremal configuration of Chen and Chen we show, by certificates
verified in exact rational arithmetic, that all $56$ such extensions have
SOS rank at most $55$: $32$ by product relations and $24$ by contraction
certificates of Xu--Yan type.  We also show, with computer-assisted
certificates and an exhaustive isomorphism search, that the extremal
configurations of $z_2(21,4)$ fall into at least six isomorphism classes;
Chen and Chen leave the structural classification of extremal augmented
graphs open.  Finally we prove a six-to-four reduction.  For irreducible
configurations whose $1$-edge graph is $C_4$-free with $m+\binom n2$ edges
it bounds the multiplicity of a shared half by $n(n-1)$ ($12$ at
$21\times4$), and it applies directly to an explicit number of the
extensions of any extremal configuration of $z_2(21,4)$, eight for that of
Chen and Chen, bringing each to at most $54$ squares, in agreement with
the certificates.  The cell count, the heredity
lemma, the rows lemma, the six-to-four identity and a $2\times2$ example
are checked in Isabelle/HOL; the remaining arguments are written proofs.
\end{abstract}

\section{Introduction}

For an oriented augmented grid $G$ on $m\times n$ cells, with a set $E_1$
of $1$-edges (cells) and a set $E_2$ of $2$-edges (pairs of cells, each
with a main half and a secondary half), let
\[
P_G=\sum_{(i,j)\in E_1}(x_iy_j)^2+\sum_{\{(i,j),(k,\ell)\}\in E_2}(x_iy_j+x_ky_\ell)^2 .
\]
We write $\SOS(P)$ for the least number of squares of real bilinear forms
whose sum is $P$, and call the display of $G$ \emph{irreducible} when
$\SOS(P_G)=|E_1|+|E_2|$, the number $t$ of displayed squares.
$\BSR(m,n)$ is the largest value of $\SOS(P)$ over all $m\times n$
biquadratic forms that are sums of squares.  In the definition of
$z_{2M}(m,n)$ \cite[Definitions~2--3]{QiYu6} the $1$-edge graph is
$C_4$-free with $z(m,n)$ edges; the $1$-edge cells, the main halves and
the distinct secondary halves form pairwise disjoint sets; $1$-edge
cells and main halves are exclusive, while a secondary half may be shared
by several $2$-edges (condition $(S_M)$).  Without sharing this is
condition $(S)$ of the second order Zarankiewicz number $z_2$
\cite{LQ}.  $H$ denotes the number of empty cells.

Qi and Yu \cite{QiYu6} ask whether $z_{2M}(21,4)=55$.  The known bounds
are $55\le z_{2M}(21,4)\le 63$.  The lower bound comes from the extremal
configuration of $z_2(21,4)=55$ of Chen and Chen \cite{CC}, which
satisfies $(S)$ and hence $(S_M)$; the upper bound is $z_{2M}\le\BSR$
together with $\BSR(21,4)\le 63$ \cite[Corollary~7]{MarinW3}.

Table~\ref{tab:new} separates what is proved here from what we use.

\begin{table}[ht]\centering
\caption{Results of this note and the known tools they rest on.}\label{tab:new}
\small
\begin{tabular}{>{\raggedright\arraybackslash}p{0.56\textwidth}>{\raggedright\arraybackslash}p{0.36\textwidth}}
\hline
\rowcolor{headbg}statement & status\\
\hline
\semrow{knownbg}heredity of irreducibility under deletion of displayed squares & known \cite[Lemma~4.17]{LQ}, \cite[Lemma~2.11]{XY}\\
\semrow{newbg}$z_{2M}(21,4)=55$ iff no irreducible $56$-square configuration (Prop.~\ref{prop:56}) & correction of \cite[Remark~10]{QiYu6}, where a deletion lowers $t$ by two\\
\semrow{knownbg}cell count $2(t-55)=1+k-H$, and $k=1$ when $t=56$, $H=0$ & known \cite[Remark~10]{QiYu6}; cf.\ \cite[eq.~(2)]{XY}\\
\semrow{newbg}irreducible $55$-square cores of an irreducible $56$-square configuration with $H=0$ (Prop.~\ref{prop:count}) & from the count and heredity; proved here\\
\semrow{knownbg}rows of an extremal $C_4$-free graph (Lemma~\ref{lem:rows}) & classical; \cite[Lemma~4.1]{CC}, \cite[Lemma~4.6]{XY}\\
\semrow{newbg}the $56$ extensions of the Chen--Chen base have $\SOS\le55$ (Thm.~\ref{thm:56}) & new (computation); methods from \cite{LQ,XY}\\
\semrow{knownbg}product-relation shortening & known \cite[Lemma~B.2]{LQ}, \cite[Lemma~2.13]{XY}\\
\semrow{knownbg}contraction criterion & known \cite[Lemma~2.15]{XY}\\
\semrow{knownbg}closure rules for irreducibility & known \cite[Definitions~4.1, 4.5; Theorem~4.3, Proposition~4.6]{LQ}, \cite[\S2]{Leb}\\
\semrow{newbg}at least six isomorphism classes of extremal configurations & new; classification open in \cite[\S5.2]{CC}\\
\semrow{newbg}six-to-four reduction for a shared half and multiplicity $\le n(n-1)$ (Lemma~\ref{lem:64}, Cor.~\ref{cor:12}) & new application of the product-relation mechanism, and its density consequence\\
\semrow{newbg}exactly $\sum_{b\in S_0}(\deg b-2)$ extensions admit a direct application of Lemma~\ref{lem:64}, each with $\SOS\le54$ (Prop.~\ref{prop:count64}) & new\\
\semrow{newbg}two cores with unique PSD Gram matrices whose fully occupied union is irreducible (Remark~\ref{rem:22}) & new example; shared-half examples without extremal $1$-edge graph in \cite[Remark~6]{QiYu6}\\
\semrow{knownbg}column bound & projection argument of \cite[Corollary~7]{MarinW3}, \cite[Theorem~9]{MarinQP}, after \cite[Corollary~2.15]{BSSV}\\
\hline
\end{tabular}
\end{table}

\section{The $56$-square reduction}

\begin{lemma}[{\cite[Lemma~4.17]{LQ}, \cite[Lemma~2.11]{XY}}]\label{lem:her}
If $P=Q+f^2$ and $\SOS(P)=t$, then $\SOS(Q)\ge t-1$.
\end{lemma}

Indeed, a representation of $Q$ with fewer than $t-1$ squares gives one of
$P$ with fewer than $t$.

\begin{proposition}\label{prop:56}
$z_{2M}(21,4)=55$ if and only if there is no irreducible
$(S_M)$-admissible configuration at $21\times4$ whose $1$-edge graph is
$C_4$-free with $27$ edges and which has exactly $56$ displayed squares.
\end{proposition}

\begin{proof}
Delete a $2$-edge of an irreducible admissible configuration with
$t\ge 57$ squares, keeping the orientation of every remaining edge.
This removes one complete square, leaves $E_1$ unchanged and keeps
$(S_M)$; by Lemma~\ref{lem:her} the result is irreducible with $t-1$
squares.  Repeating, an irreducible admissible configuration of any
length $t\ge 56$ yields one of length $56$.  Since an irreducible
admissible $55$-square configuration exists \cite{CC}, the statement
follows.  This corrects the deletion count in \cite[Remark~10]{QiYu6}:
a $2$-edge contributes one displayed square, so there is no separate
case $t=57$.
\end{proof}

The same descent to the next size is used for $z_2$ in
\cite[proof of Corollary~4.19]{LQ} and \cite[proof of Theorem~4.4]{XY}.
The cell count below is \cite[Remark~10]{QiYu6} (compare
\cite[eq.~(2)]{XY} for $z_2$); we include its short proof because the
rest of the note depends on it.

\begin{proposition}\label{prop:count}
Let $G$ be admissible at $21\times4$ with $|E_1|=27$, $e_2=|E_2|$ two-edges,
$s$ distinct secondary halves, $k=e_2-s$ and $t=27+e_2$.  Then
\[
2(t-55)=1+k-H ,
\]
and if $t=56$ and $H=0$, then $k=1$ \cite[Remark~10]{QiYu6}.  If moreover
$G$ is irreducible, deleting either of the two
$2$-edges at the shared secondary half leaves an irreducible configuration
satisfying $(S)$ with $55$ squares and one hole, that is, an extremal
configuration of $z_2(21,4)=55$.
\end{proposition}

\begin{proof}
The $1$-edge cells, the main halves and the distinct secondary halves are
pairwise disjoint, of sizes $27$, $e_2$ and $s$, so $84=27+e_2+s+H$.
Substituting $s=e_2-k$ and $e_2=t-27$ gives the identity.  At $t=56$,
$H=0$ it reads $k=1$: exactly one secondary cell is used by two
$2$-edges and every other by one.  Deleting either of these two edges
makes its exclusive main half the unique hole; every remaining secondary
cell has multiplicity one, so $(S)$ holds, and irreducibility follows from
Lemma~\ref{lem:her}.
\end{proof}

So every irreducible admissible $56$-square configuration with $H=0$ is an
extremal $55$-square configuration plus one $2$-edge joining its hole to a
half of one of its $2$-edges.  For $H\ge1$ one has $k=H+1\ge2$, and
deleting one $2$-edge can leave sharing.

\section{The $56$ extensions of the Chen--Chen base}

Let $G_0$ be the $21$-row configuration obtained in \cite[\S4.2]{CC} by
the deletion of \cite[Proposition~4.10]{CC}, with rows renumbered
consecutively.  It has $27$ one-edges, $28$ two-edges and the unique hole
$p=(10,1)$.  For each half $c$ of an old $2$-edge, orient that edge with
secondary half $c$ and adjoin the $2$-edge with main half $p$ and
secondary half $c$.  This gives an $(S_M)$-admissible display with $56$
squares and $H=0$; orientation does not change $P_G$, and since a
$1$-edge cell cannot be shared these $56$ displays are all the extensions
of $G_0$ by one $2$-edge at $p$.

\begin{theorem}\label{thm:56}
For each of the $56$ extensions, $\SOS(P_G)\le 55$.
\end{theorem}

\begin{proof}
Write $f_1,\dots,f_{56}$ for the displayed forms and $z$ for the vector of
the $84$ cell monomials, $f=Cz$.

\emph{(i) Product relations, $32$ extensions.}  The certificate gives
rational numbers $c_{rs}$, not all zero, with
$\sum_{1\le r<s\le 56}c_{rs}f_rf_s=0$ identically.  By the product-relation
lemma \cite[Lemma~B.2]{LQ} (see also \cite[Lemma~2.13]{XY}), $P_G$ is a
sum of at most $55$ squares: with the symmetric matrix $K$,
$K_{rs}=K_{sr}=c_{rs}/2$, $K_{rr}=0$, and $\lambda=\lambda_{\max}(K)>0$,
the matrix $I-K/\lambda$ is positive semidefinite, singular, and represents
$P_G$ in the displayed basis.

\emph{(ii) Contraction certificates, $24$ extensions,} among them
$c=(17,3)$.  Here the $1596$ products $f_rf_s$, $1\le r\le s\le56$, are
linearly independent.  Let $W$ be the span of the rows of $C$; each
displayed form has an exclusive cell, so $C$ has rank $56$, and a Gram
matrix with range in $W$ has the form $C^TAC$.  Product independence forces
$A=I$.  Hence the displayed Gram matrix is an extreme point of the Gram
spectrahedron of $P_G$ and the only feasible Gram matrix with range in
$W$, so a representation with fewer than $56$ squares uses at least one
form outside the displayed span (compare \cite[Proposition~5.4]{XY}).
The certificate \cite[Lemma~2.15]{XY}
specifies a rational affine family $Q(h)=Q_0+\sum_j h_jB_j$ and a rational
matrix $L$ of full column rank with $z^TLQ(h)L^Tz\equiv P_G$ for every
$h$: the checker verifies the identity at $h=0$ and that every direction
$B_j$ has zero polynomial image.  After the recorded permutation,
\[
Q(h)=\begin{pmatrix}A(h)&D(h)\\ D(h)^T&E(h)\end{pmatrix},
\]
with $A(h)$ of order $55$.  The exact bounds of \cite[Lemma~2.15]{XY}
prove $A(h)\succ0$ on the box $\|h\|_\infty\le\varepsilon=10^{-10}$ and a
zero $h_*$ of the Schur complement $E-D^TA^{-1}D$ in that box.  Then
$LQ(h_*)L^T$ is positive semidefinite of rank exactly $55$, and
$\SOS(P_G)\le 55$.
\end{proof}

\begin{figure}[ht]\centering
\includegraphics[width=\textwidth]{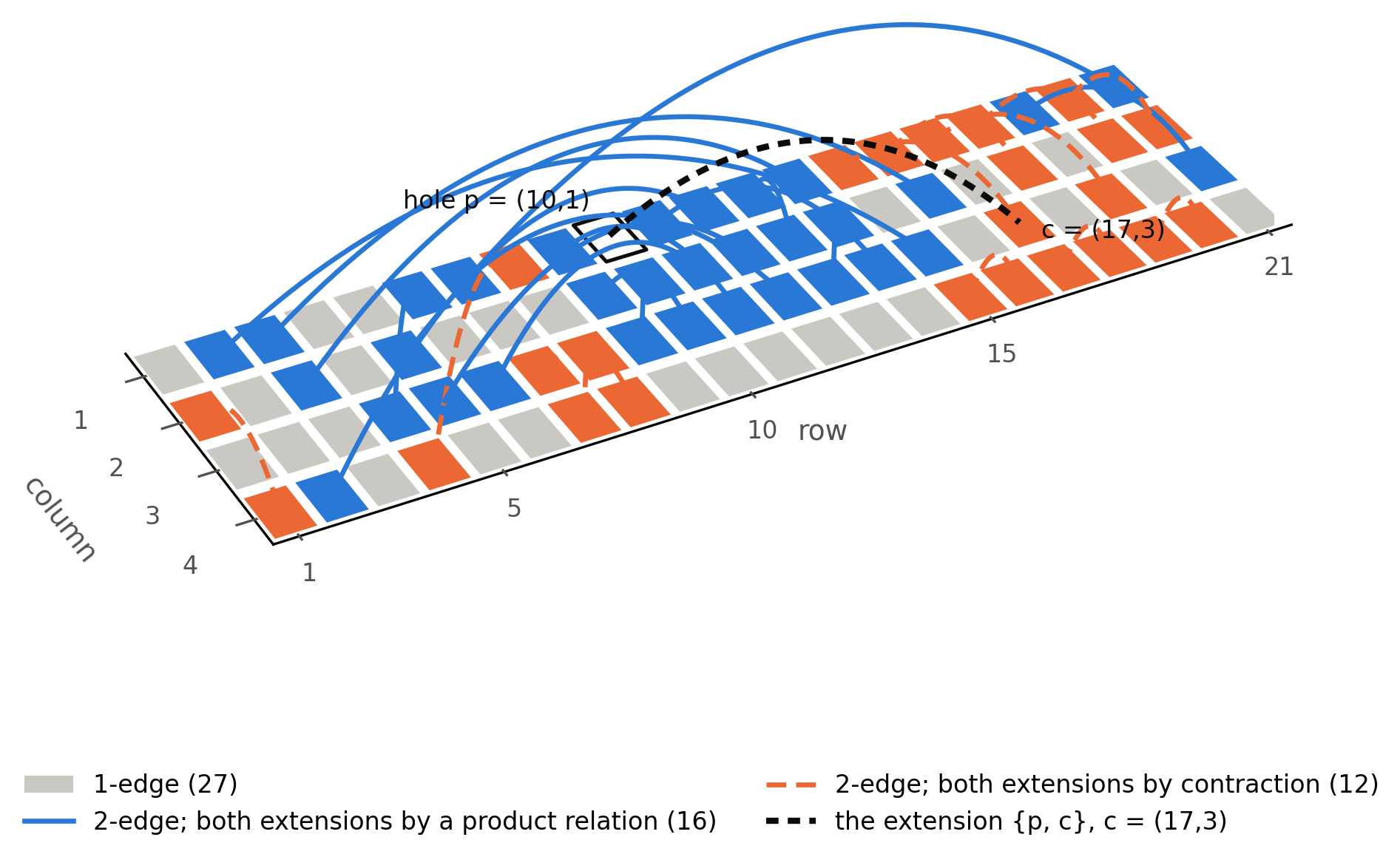}
\caption{The base $G_0$ on the $21\times4$ grid.  Grey tiles: $1$-edges.
Each $2$-edge is an arc between its two cells, solid blue if both of its
extensions are reduced by a product relation and dashed orange if both are
reduced by a contraction certificate.  The black dashed arc is the
extension from the hole $p=(10,1)$ to $c=(17,3)$.  The complete pairings
are in \texttt{chen\_chen\_base55.json} of the supplement, and the plotting
script is included.}\label{fig:3d}
\end{figure}

Table~\ref{tab:ext} lists the certificate of each extension and
Figure~\ref{fig:3d} shows where the two kinds lie on the grid.  The type
depends only on the old $2$-edge: both halves of an old edge give the same
kind of certificate ($16$ edges by product relations, $12$ by
contraction).

\begin{remark}[Scope of the certificates]
All acceptance checks use exact rational arithmetic implemented with
arbitrary-precision integers; the verifier rebuilds each target
polynomial from the base and rejects three deliberately altered
certificates.  The rank-$55$ Gram matrix of item (ii) is proved to exist;
its entries are not asserted to be rational, and it is not written down.
The certificates give upper bounds only, for the $56$ extensions of
$G_0$ and for configurations obtained from them by relabelling rows and
columns.  They do not give the minimum number of squares.  Apart from one
extension of each configuration $N_c$ of Section~\ref{sec:other} (since
$G_0+\{p,c\}=N_c+\{d,c\}$), they say nothing about extensions of other
extremal configurations.
\end{remark}

\begin{longtable}{lllrrr}
\caption{The $56$ extensions $\{p,c\}$, $p=(10,1)$, grouped by the old $2$-edge $\{c,c^{\prime}\}$ of the base. For a product relation: the number of products $f_rf_s$ in the certified relation. For a contraction: the order $N_Q$ of the reduced Gram block $Q(h)$, the contraction constant $\theta$ and the bound $\alpha+\theta\varepsilon$ checked against $\varepsilon=10^{-10}$. Decimal entries are rounded for display; the acceptance comparisons use the exact rational values in the supplementary data.}\label{tab:ext}\\
\hline
\rowcolor{headbg}$2$-edge & $c$ & certificate & $N_Q$ & $\theta$ & $\alpha+\theta\varepsilon$\\
\hline
\endfirsthead
\multicolumn{6}{l}{\small Table~\ref{tab:ext}, continued}\\
\hline
\rowcolor{headbg}$2$-edge & $c$ & certificate & $N_Q$ & $\theta$ & $\alpha+\theta\varepsilon$\\
\hline
\endhead
\semrow{conbg}$e_{1}$ & (1,2) & contraction & 76 & 0.235 & $2.4\cdot 10^{-11}$\\
\semrow{conbg} & (1,4) & contraction & 76 & 0.040 & $5.0\cdot 10^{-12}$\\
\semrow{relbg}$e_{2}$ & (2,1) & product relation (15 products) &  &  & \\
\semrow{relbg} & (15,3) & product relation (13 products) &  &  & \\
\semrow{relbg}$e_{3}$ & (2,4) & product relation (4 products) &  &  & \\
\semrow{relbg} & (12,1) & product relation (3 products) &  &  & \\
\semrow{relbg}$e_{4}$ & (3,1) & product relation (3 products) &  &  & \\
\semrow{relbg} & (16,2) & product relation (4 products) &  &  & \\
\hline
\semrow{relbg}$e_{5}$ & (3,2) & product relation (4 products) &  &  & \\
\semrow{relbg} & (14,1) & product relation (3 products) &  &  & \\
\semrow{relbg}$e_{6}$ & (4,3) & product relation (4 products) &  &  & \\
\semrow{relbg} & (6,1) & product relation (3 products) &  &  & \\
\semrow{conbg}$e_{7}$ & (4,4) & contraction & 75 & 0.053 & $5.3\cdot 10^{-12}$\\
\semrow{conbg} & (8,1) & contraction & 75 & 0.277 & $2.8\cdot 10^{-11}$\\
\semrow{relbg}$e_{8}$ & (5,2) & product relation (4 products) &  &  & \\
\semrow{relbg} & (21,1) & product relation (3 products) &  &  & \\
\hline
\semrow{relbg}$e_{9}$ & (5,3) & product relation (6 products) &  &  & \\
\semrow{relbg} & (12,2) & product relation (7 products) &  &  & \\
\semrow{relbg}$e_{10}$ & (6,3) & product relation (4 products) &  &  & \\
\semrow{relbg} & (11,2) & product relation (5 products) &  &  & \\
\semrow{relbg}$e_{11}$ & (7,1) & product relation (10 products) &  &  & \\
\semrow{relbg} & (12,3) & product relation (11 products) &  &  & \\
\semrow{conbg}$e_{12}$ & (7,3) & contraction & 77 & 0.009 & $1.2\cdot 10^{-12}$\\
\semrow{conbg} & (8,4) & contraction & 77 & 0.061 & $6.3\cdot 10^{-12}$\\
\hline
\semrow{conbg}$e_{13}$ & (7,4) & contraction & 77 & 0.049 & $5.0\cdot 10^{-12}$\\
\semrow{conbg} & (8,3) & contraction & 77 & 0.103 & $1.1\cdot 10^{-11}$\\
\semrow{relbg}$e_{14}$ & (9,1) & product relation (3 products) &  &  & \\
\semrow{relbg} & (11,3) & product relation (4 products) &  &  & \\
\semrow{relbg}$e_{15}$ & (9,2) & product relation (10 products) &  &  & \\
\semrow{relbg} & (10,3) & product relation (10 products) &  &  & \\
\semrow{relbg}$e_{16}$ & (9,3) & product relation (11 products) &  &  & \\
\semrow{relbg} & (10,2) & product relation (2 products) &  &  & \\
\hline
\semrow{relbg}$e_{17}$ & (11,1) & product relation (3 products) &  &  & \\
\semrow{relbg} & (14,3) & product relation (4 products) &  &  & \\
\semrow{relbg}$e_{18}$ & (13,1) & product relation (3 products) &  &  & \\
\semrow{relbg} & (13,2) & product relation (4 products) &  &  & \\
\semrow{relbg}$e_{19}$ & (13,3) & product relation (9 products) &  &  & \\
\semrow{relbg} & (14,2) & product relation (10 products) &  &  & \\
\semrow{conbg}$e_{20}$ & (15,1) & contraction & 75 & 0.020 & $2.3\cdot 10^{-12}$\\
\semrow{conbg} & (17,3) & contraction & 75 & 0.026 & $3.3\cdot 10^{-12}$\\
\hline
\semrow{conbg}$e_{21}$ & (15,4) & contraction & 76 & 0.207 & $2.1\cdot 10^{-11}$\\
\semrow{conbg} & (16,4) & contraction & 76 & 0.090 & $9.2\cdot 10^{-12}$\\
\semrow{conbg}$e_{22}$ & (16,1) & contraction & 75 & 0.018 & $1.8\cdot 10^{-12}$\\
\semrow{conbg} & (18,2) & contraction & 75 & 0.046 & $4.6\cdot 10^{-12}$\\
\semrow{conbg}$e_{23}$ & (17,1) & contraction & 75 & 0.027 & $2.7\cdot 10^{-12}$\\
\semrow{conbg} & (19,3) & contraction & 75 & 0.039 & $4.0\cdot 10^{-12}$\\
\semrow{conbg}$e_{24}$ & (17,4) & contraction & 76 & 0.106 & $1.1\cdot 10^{-11}$\\
\semrow{conbg} & (18,4) & contraction & 76 & 0.118 & $1.3\cdot 10^{-11}$\\
\hline
\semrow{conbg}$e_{25}$ & (18,1) & contraction & 75 & 0.028 & $2.9\cdot 10^{-12}$\\
\semrow{conbg} & (20,2) & contraction & 75 & 0.053 & $5.7\cdot 10^{-12}$\\
\semrow{relbg}$e_{26}$ & (19,1) & product relation (10 products) &  &  & \\
\semrow{relbg} & (21,3) & product relation (9 products) &  &  & \\
\semrow{conbg}$e_{27}$ & (19,4) & contraction & 76 & 0.183 & $1.9\cdot 10^{-11}$\\
\semrow{conbg} & (20,4) & contraction & 76 & 0.066 & $6.8\cdot 10^{-12}$\\
\semrow{conbg}$e_{28}$ & (20,1) & contraction & 75 & 0.027 & $2.7\cdot 10^{-12}$\\
\semrow{conbg} & (21,2) & contraction & 75 & 0.049 & $5.2\cdot 10^{-12}$\\
\hline
\end{longtable}

\section{Other extremal configurations}\label{sec:other}

Replacing an old $2$-edge $\{d,c\}$ of $G_0$ by $\{p,c\}$ gives a
$55$-square configuration $N_c$ with its hole at $d$.  Of the $56$
configurations $N_c$, sixteen have a product relation among their
displayed forms, hence are reducible; the closure below decides six as
irreducible; the other $34$ are undecided.  For the six choices of $c$ in
Table~\ref{tab:neigh}, the closure
rules of \cite[Definitions~4.1 and~4.5]{LQ} and \cite[\S2]{Leb}, sound by
\cite[Theorem~4.3, Proposition~4.6]{LQ} (coefficient
identities on rectangles, the empty cell, Cauchy--Schwarz saturation and
substitution of equal vectors) fix all $\binom{83}{2}+83=3486$ products
among the $83$ occupied cell vectors, diagonal included, in every positive
semidefinite Gram matrix (the products with the empty cell vanish), and
the resulting Gram matrix has rank $55$.  So each of the six is
irreducible, hence extremal for $z_2(21,4)$.  Each certificate is accepted
by a separate checker that does not import the generating program.

\begin{table}[ht]\centering
\caption{The six irreducible $55$-square configurations obtained from the base by replacing the old $2$-edge $\{d,c\}$ by $\{p,c\}$. Class: isomorphism class under row and column permutations ($0$ is the class of the base).}\label{tab:neigh}
\begin{tabular}{llcc}
\hline
\rowcolor{headbg}kept half $c$ & new hole $d$ & products fixed among occupied cells & class\\
\hline
\semrow{newbg}(9,3) & (10,2) & all 3486 $=\binom{83}{2}+83$ & 1\\
\semrow{knownbg}(11,3) & (9,1) & all 3486 $=\binom{83}{2}+83$ & 0\\
\semrow{newbg}(9,2) & (10,3) & all 3486 $=\binom{83}{2}+83$ & 2\\
\semrow{newbg}(10,3) & (9,2) & all 3486 $=\binom{83}{2}+83$ & 3\\
\semrow{newbg}(13,3) & (14,2) & all 3486 $=\binom{83}{2}+83$ & 4\\
\semrow{newbg}(13,2) & (13,1) & all 3486 $=\binom{83}{2}+83$ & 5\\
\hline
\end{tabular}
\end{table}

The configuration with $c=(11,3)$ is $G_0$ with rows $9$ and $10$
exchanged.  The other five are pairwise non-isomorphic and not
isomorphic to $G_0$ (exhaustive search over the $24$ column permutations
and all row bijections).  A short invariant already separates one of
them: $G_0$ has two $2$-edges with both halves in one row, $\{(1,2),(1,4)\}$
and $\{(13,1),(13,2)\}$, while the configuration with $c=(13,2)$ has only
one.  Hence the extremal configurations of $z_2(21,4)$ fall into at least
six isomorphism classes, and the case $H=0$ is a question about all of
them.

\section{Structural tools}

\subsection{Rows of an extremal $1$-edge graph}

\begin{lemma}[{\cite[Lemma~4.1]{CC}, \cite[Lemma~4.6]{XY}}]\label{lem:rows}
Let $E\subseteq[m]\times[n]$ be $C_4$-free with $|E|=m+\binom n2$.  Then
every row has degree $1$ or $2$, and every pair of columns lies in exactly
one row.  In particular this holds for every $1$-edge graph of
$z_{2M}(21,4)$, where $|E_1|=27=21+\binom42$.
\end{lemma}

\begin{proof}
Let $d_i$ be the degree of row $i$.  Since $E$ is $C_4$-free, distinct rows
contain distinct pairs of columns, so $\sum_i\binom{d_i}{2}\le\binom n2$.
With $d\le 1+\binom d2$ for every $d\ge0$,
\[
m+\tbinom n2=\sum_i d_i\le\sum_i\Bigl(1+\tbinom{d_i}{2}\Bigr)\le m+\tbinom n2 .
\]
Equality throughout gives $d_i=1+\binom{d_i}2$, that is $d_i\in\{1,2\}$,
and $\sum_i\binom{d_i}2=\binom n2$, so the pairs of columns in the rows
are all of them, each once.
\end{proof}

This is the classical equality case of the counting behind
$z(m,n)\le m+\binom n2$; we include the proof because it is short and
checked in Isabelle (Section~\ref{sec:formal}).

\subsection{The six-to-four reduction}

\begin{lemma}\label{lem:64}
Suppose an $(S_M)$-admissible display of length $t$ contains the
$2$-edges $f_i=x_iy_c+q$ and $f_j=x_jy_c+q$ with the same secondary
monomial $q$ and $i\ne j$, and, for some column $b\ne c$ and row
$k\notin\{i,j\}$, the $1$-edges $a_i=x_iy_b$, $a_j=x_jy_b$, $u=x_ky_c$,
$v=x_ky_b$.  Then $\SOS(P_G)\le t-2$.
\end{lemma}

\begin{proof}
Put $F_\pm=(f_i\pm f_j)/\sqrt2$ and $A_\pm=(a_i\pm a_j)/\sqrt2$.  Then
$F_-=(x_i-x_j)y_c/\sqrt2$ and $A_-=(x_i-x_j)y_b/\sqrt2$, so $F_-v=A_-u$,
and
\[
f_i^2+f_j^2+a_i^2+a_j^2+u^2+v^2=F_+^2+A_+^2+(F_-+v)^2+(A_--u)^2 .
\]
Replace the six summands by the four on the right and keep the others.
\end{proof}

\begin{corollary}\label{cor:12}
Let $G$ be an irreducible $(S_M)$-admissible configuration on $m\times n$
cells whose $1$-edge graph is $C_4$-free with $|E_1|=m+\binom n2$.  Then a
shared secondary half has at most $n-1$ main halves in any one column,
hence multiplicity at most $n(n-1)$.  At $21\times4$, with $|E_1|=27$, the
multiplicity is at most $12$.
\end{corollary}

\begin{proof}
Fix a column $c$ and the rows $i_1,i_2,\dots$ of the main halves in $c$ of
the $2$-edges at one secondary half.  The $1$-edges of row $i_r$ avoid
column $c$, which holds a main half, and by Lemma~\ref{lem:rows} each row
has one or two $1$-edges.  If rows $i_r\ne i_s$ had $1$-edges in a common
column $b$, then by Lemma~\ref{lem:rows} the pair $\{b,c\}$ lies in a row
$k$ with $1$-edges at $(k,b)$ and $(k,c)$; $k\notin\{i_r,i_s\}$ because
$(i_r,c)$ and $(i_s,c)$ are main halves, and Lemma~\ref{lem:64} would
reduce the display.  So the column sets of $1$-edges of these rows are
nonempty and pairwise disjoint subsets of the $n-1$ other columns.
\end{proof}

\begin{remark}\label{rem:transpose}
The transposed form of Lemma~\ref{lem:64}, with two main halves $(i,c)$,
$(i,c')$ in one row, needs a row $b$ with $1$-edges at $(b,c)$, $(b,c')$
and $(b,k)$ for a third column $k$.  Under the hypothesis of
Corollary~\ref{cor:12} every row has at most two $1$-edges
(Lemma~\ref{lem:rows}), so the transposed form never applies there.
\end{remark}

\begin{proposition}\label{prop:count64}
Let $G$ be an extremal $55$-square configuration of $z_2(21,4)$ with its
unique hole $p=(r_0,c_0)$, and let $S_0$ be the set of columns of the
$1$-edges of row $r_0$.  Among the $56$ extensions of $G$ by one $2$-edge
at $p$, exactly
\[
\sum_{b\in S_0}\bigl(\deg b-2\bigr)
\]
admit a direct application of Lemma~\ref{lem:64} to the displayed
representation, where $\deg b$ is the number of $1$-edges in column $b$.
Each of these extensions has SOS rank at most $54$.
\end{proposition}

\begin{proof}
Since $G$ satisfies $(S)$, the only shared cell in the extension $\{p,c\}$
is $c$, common to the new $2$-edge and the old one $\{d,c\}$; condition
$(S_M)$ forces $c$ to be secondary in both, with main halves $p$ and $d$.
The extension keeps $E_1$, so its rows still have at most two $1$-edges,
and by Remark~\ref{rem:transpose} Lemma~\ref{lem:64} applies only
when $d=(r,c_0)$ lies in column $c_0$ and rows $r_0$ and $r$ have $1$-edges
in a common column $b$; then the row $k_b$ of the pair $\{b,c_0\}$
(Lemma~\ref{lem:rows}) supplies $(k_b,b)$ and $(k_b,c_0)$, and
$k_b\notin\{r_0,r\}$ since $(r_0,c_0)$ is the hole and $(r,c_0)$ a half.
Conversely, let $r\ne r_0$ have a $1$-edge in some $b\in S_0$ and
$r\ne k_b$ for all $b$.  Then $(r,c_0)$ is not a $1$-edge (otherwise $r$
would contain the pair $\{b,c_0\}$), and not a hole, so it is a half $d$,
and its extension is reduced.  Each such row gives one extension.  The rows
with a $1$-edge in some $b\in S_0$ number $\sum_{b\in S_0}\deg b-(|S_0|-1)$,
since a row with two columns of $S_0$ is the row $r_0$ of that pair;
removing $r_0$ and the $|S_0|$ distinct rows $k_b$ leaves
$\sum_{b\in S_0}(\deg b-2)$.  Lemma~\ref{lem:64} replaces six of the $56$
displayed squares by four, so each counted extension has at most $54$
squares; this is an upper bound, not the minimum.
\end{proof}

For $G_0$, $S_0=\{4\}$ and $\deg 4=10$, so Lemma~\ref{lem:64} reduces $8$
extensions by hand: $c=(2,4)$, $(16,2)$, $(3,2)$, $(4,3)$, $(5,2)$, $(11,3)$,
$(14,3)$, $(13,2)$, all with $b=4$ and $k_b=5$.  This agrees with
Table~\ref{tab:ext}: all eight are certified there by product relations
with four products each, the shape of the relation $F_-v-A_-u=0$ of
Lemma~\ref{lem:64}.  The six configurations of Table~\ref{tab:neigh} share
the $1$-edge graph of $G_0$, and in each of them the row of the new hole
again has $S_0=\{4\}$, so the count is $\deg 4-2=8$ in every case.

\subsection{The column bound}

For an SOS biquadratic form $P$ write $P(x,e_j)=x^TQ_jx$.  Then
\[
\SOS(P)\le\operatorname{rank}Q_j+\BSR(m,n-1).
\]
This is the rank-sensitive form, for a single column, of the
coefficient-vector projection argument used in \cite[proof of
Corollary~7]{MarinW3} and \cite[proof of Theorem~9]{MarinQP}; the earlier
dimension-only bound is a case of the parallel-line covering bound of
\cite[Corollary~2.15]{BSSV}.  In an SOS representation
$P=\|\sum_{i,\ell}v_{i\ell}x_iy_\ell\|^2$,
the vectors $v_{ij}$, $1\le i\le m$, span a space of dimension
$\operatorname{rank}Q_j$, and the orthogonal projection onto it and onto
its complement splits $P$ into a sum of at most $\operatorname{rank}Q_j$
squares and an SOS form free of $y_j$.  If $(i_0,j)$ is an empty cell,
$v_{i_0j}=0$ in every representation (its squared norm is the zero
coefficient of $x_{i_0}^2y_j^2$), so $\operatorname{rank}Q_j\le 20$ and,
with $\BSR(21,3)=42$ \cite[Theorem~1]{MarinW3}, $\SOS(P_G)\le 62$.  This
does not exclude an irreducible $56$-square display, since $56\le62$; and
for $H=0$ it gives only $21+42=63$.

\begin{remark}[A $2\times2$ example]\label{rem:22}
Let $a=x_1y_1$, $b=x_1y_2$, $c=x_2y_1$, $d=x_2y_2$.  The forms
$a^2+(b+c)^2$ and $a^2+(b+d)^2$ each have a unique positive semidefinite
Gram matrix (the zero diagonal entry of the empty cell fixes the only
free parameter), and their union $a^2+(b+c)^2+(b+d)^2$, with every cell
occupied, is not a sum of two squares: its matrix in $y$ has determinant
$2x_1^4+2x_1^3x_2+2x_1^2x_2^2+2x_1x_2^3+x_2^4$, which is not the square of
a quadratic form.  So full occupancy and two cores with unique Gram
matrices do not force reducibility of their union; a proof for the
$21\times4$ class must use further hypotheses, the extremal density of the
$1$-edge graph being one possible source.  Irreducible configurations with
a shared half but with an empty $1$-edge graph are given in
\cite[Remark~6]{QiYu6}.
\end{remark}

\section{Numerical evidence outside Theorem~\ref{thm:56}}

We also examined two $56$-square configurations with one empty cell, with
the main halves of each shared star in one column.  The first has
$(e_1,e_2,s,k)=(27,29,27,2)$; the second has $(26,30,27,3)$ and lies in
the relaxed class with fewer than $27$ one-edges.  Our recorded
floating-point searches report approximate $55$-square decompositions,
with coefficient residuals about $9.29\cdot10^{-12}$ and
$5.35\cdot10^{-12}$.  The supplementary files keep summary records only,
not the complete targets or computed factors needed to reproduce these
residuals.  These observations are numerical evidence, not certified
upper bounds, and are not used in Theorem~\ref{thm:56}.

\section{Formal verification}\label{sec:formal}

The Isabelle/HOL session \texttt{Z2M\_Hole\_Extensions}, folder
\href{https://github.com/karlesmarin/formalizations/tree/390bf7d7247a3320f4fcc04d1c6cc606e877eb53/Isabelle/Z2MHoleExtensions}{\texttt{Isabelle/Z2MHoleExtensions}}
of \url{https://github.com/karlesmarin/formalizations}, commit \texttt{390bf7d} (a self-contained
copy is in the supplement), formalizes the elementary statements listed below (built without \texttt{sorry} and
without \texttt{quick\_and\_dirty}; a copy with a false lemma appended
fails, as it should; the supplement copy is the session
\texttt{MQY\_Pieces}, with its build log).  The heredity theorem is
stated for arbitrary families of real-valued functions.  The passage from
admissible configurations to the premises of the counting and heredity
lemmas, the graph-theoretic application of the six-to-four identity,
Proposition~\ref{prop:56}, the core statement of Proposition~\ref{prop:count},
Corollary~\ref{cor:12}, Remark~\ref{rem:transpose}, Proposition~\ref{prop:count64}, the column
bound and the Gram uniqueness in Remark~\ref{rem:22} remain in the written proofs.

\begin{table}[ht]\centering
\small
\begin{tabular}{>{\raggedright\arraybackslash}p{0.4\textwidth}>{\raggedright\arraybackslash}p{0.52\textwidth}}
\hline
\rowcolor{headbg}statement & Isabelle theorem\\
\hline
\rowcolor{white}cell identity of Proposition~\ref{prop:count} & \texttt{display\_cell\_count}, \texttt{cell\_count\_21x4}\\
\rowcolor{stripebg}$t=56$, $H=0$ gives $k=1$ & \texttt{no\_hole\_one\_share}\\
\rowcolor{white}heredity, Lemma~\ref{lem:her} & \texttt{irreducible\_delete\_square}\\
\rowcolor{stripebg}Lemma~\ref{lem:rows} and the case $21\times4$ & \texttt{c4\_free\_extremal\_rows}, \texttt{every\_column\_pair\_once\_21x4}\\
\rowcolor{white}the identity of Lemma~\ref{lem:64} & \texttt{six\_to\_four}\\
\rowcolor{stripebg}the $2\times2$ example (not two squares; three) & \texttt{ex22\_not\_two\_squares}, \texttt{ex22\_three\_squares}\\
\hline
\end{tabular}
\end{table}

The certificates of Theorem~\ref{thm:56} and Table~\ref{tab:neigh} are
checked by the exact verifiers of the supplementary data, not in
Isabelle.

\section{Data and reproducibility}

The supplementary data contain the base $G_0$, the $56$ extensions, the
$32$ product-relation certificates, the $24$ contraction certificates,
the six closure certificates of Table~\ref{tab:neigh}, the exact
verification programs, the Isabelle session, and a README with the
commands, the expected acceptance counts and the checksums of all files.
The verification programs use only the Python standard library.

\section{What remains open}

Theorem~\ref{thm:56} rules out the $56$ extensions of $G_0$ without empty
cells.  By Proposition~\ref{prop:count}, an irreducible admissible
$56$-square configuration with $H=0$ is obtained by adjoining one
$2$-edge to an extremal $55$-square configuration, and neither of its two
cores is isomorphic to $G_0$.  Excluding all of them requires a result
for the other extremal configurations, of which there are at least five
more classes.  For $H\ge1$ deleting a $2$-edge can leave sharing, and
those cases remain open.

\section*{Acknowledgment}
The first author used an AI assistant during the computational
exploration and drafting.  The exact certificates are checked by the
verification programs of the supplementary data.


\begin{thebibliography}{9}

\bibitem{CC} Z.~Chen, Y.~Chen,
Recursive-Line Zarankiewicz Numbers with Four Columns,
arXiv:2609.11093v2, 2026.

\bibitem{BSSV} G.~Blekherman, R.~Sinn, G.~G.~Smith, M.~Velasco,
Sums of squares and quadratic persistence on real projective varieties,
J. Eur. Math. Soc. 24 (2022), 925--965; statement numbers as in
arXiv:1902.02754v2.

\bibitem{Leb} N.~Lebedev,
Asymptotic equivalence and exact values for second-order Zarankiewicz
numbers, arXiv:2610.05442v1, 2026.

\bibitem{LQ} J.~L\"ofberg, L.~Qi,
Second Order Zarankiewicz Number, arXiv:2608.30555v3, 2026.

\bibitem{MarinW3} C.~Mar\'in,
For $m\ge2$, every $m\times3$ biquadratic sum of squares is a sum of $2m$
squares, version~2, Zenodo, 2026, doi:10.5281/zenodo.24051362.

\bibitem{MarinQP} C.~Mar\'in,
Quadratic persistence of Segre varieties and a linear gap between
biquadratic SOS rank and Zarankiewicz bounds, Zenodo, 2026,
doi:10.5281/zenodo.24036320.

\bibitem{QiYu6} L.~Qi, G.~Yu,
A Second Order Extended Zarankiewicz Number for Biquadratic SOS Rank,
version~6, Zenodo, 2026, doi:10.5281/zenodo.24156935.

\bibitem{XY} Y.~Xu, X.~Yan,
Exact Values, Extremal Classifications, and Sum-of-Squares Reductions for
Second-Order Zarankiewicz Numbers, arXiv:2609.18429v1, 2026.

\end{thebibliography}
\end{document}